\documentclass[11pt, reqno]{amsart}

\usepackage[hmargin=0.8in,twosideshift=0in,height=8.6in]{geometry}
\usepackage{amssymb,amsthm}
\usepackage{delarray,verbatim}
\usepackage{natbib}

\usepackage{ifpdf}
\ifpdf
\usepackage[pdftex]{graphicx}
\else
\usepackage[dvips]{graphicx}
\fi

\usepackage{ifpdf}
\usepackage{color}
\definecolor{webgreen}{rgb}{0,.5,0}
\definecolor{webbrown}{rgb}{.8,0,0}
\definecolor{emphcolor}{rgb}{0.95,0.95,0.95}

\usepackage{hyperref}
\hypersetup{%
          colorlinks=true,
          linkcolor=webbrown,
          filecolor=webbrown,
          citecolor=webgreen,
          breaklinks=true}
\ifpdf
\hypersetup{pdftex,
            pdfstartview=FitH, 
            bookmarksopen=true,
            bookmarksnumbered=true
}
\else
\hypersetup{dvips}
\fi

\usepackage{multicol}
\usepackage{multirow}

\newtheorem{thm}{Theorem}[section]
\newtheorem{lemma}{Lemma}[section]
\newtheorem{remark}{Remark}[section]
\newtheorem{cor}{Corollary}[section]

\numberwithin{equation}{section}

\newcommand{\F}{\mathcal{F}}
\newcommand{\C}{\mathcal{C}}

\renewcommand{\P}{\mathbb{P}}
\newcommand{\E}{\mathbb{E}^x}
\newcommand{\Eb}{\mathbb{E}}

\newcommand{\eps}{\varepsilon}
\newcommand{\R}{\mathbb{R}}
\renewcommand{\S}{\mathcal{S}}

\renewcommand{\alpha}{r}

\newcommand{\A}{\mathcal{A}}

\title[On the Finite Time Horizon American Options for Jump Diffusions]{A Proof of the Smoothness of the Finite Time 
Horizon American Put Option for Jump Diffusions}

\author{Erhan Bayraktar }
\address[E. Bayraktar]{Department of
  Mathematics, University of Michigan, Ann Arbor, MI 48109}
\email{erhan@umich.edu}
\thanks{This research  is supported in part by the National Science Foundation. }
\thanks{I would like to thank the two referees and the anonymous AE for their insightful comments.
 I also would like to express my gratitude to Farid Aitsahlia, Tom Bielecki, Masahiko Egami, Sebastian Jaimungal, Sergei Levendorskii, Mihai Sirbu, Mete Soner,  Hao Xing, Virginia R. Young, Thaleia Zariphopoulou and Gordan Zitkovic for their feedback.}


\begin{document}

\begin{abstract}
We give a new proof of the fact that the value function of the
finite time horizon American put option for a jump diffusion, when
the jumps are from a compound Poisson process, is the classical
solution of a free boundary equation. We also show that the value function is $C^1$ across
the optimal stopping boundary. Our proof, which only uses the classical
theory of parabolic partial differential equations of
\cite{friedman-parabolic,friedmansde}, is an alternative to the proof that uses 
the theory of vicosity solutions (see \cite{pham97}). This new proof relies on
constructing a monotonous sequence of functions, each of which is a value
function of an optimal stopping problem for a geometric Brownian motion, converging to the value function of the American put option for the jump diffusion uniformly and exponentially fast. This sequence is constructed by iterating a functional operator that
maps a certain class of convex functions to classical solutions of corresponding free boundary equations. On the other hand, since the
approximating sequence converges to the value function
exponentially fast, it naturally leads to a good numerical
scheme. 
\end{abstract}
\maketitle
\section{Introduction}
Let $(\Omega,\F,\P)$ be a complete probability space hosting a
Wiener process $W=\{W_t;t \geq 0\}$ and a Poisson random measure $N$
on $\R_+ \times \R_+$, with mean measure $\lambda \nu(dx) dt$ (in
which $\nu$ is a probability measure on $\R_+$), independent of the
Wiener process. We will consider a Markov process $S=\{S_t; t \geq
0\}$ of the form
\begin{equation}\label{eq:dyn}
dS_t=\mu S_t dt+\sigma S_t dW_t+S_{t-} \int_{\R_+} (z-1) N(dt,dz).
\end{equation}
 In this model, if the stock price jumps at time $t$, then it moves from $S_{t-}$ to $S_t=Z  S_{t-}$, in which $Z$ is a positive random variable whose distribution is given by $\nu$.  Note that when $Z <1$ the stock price  jumps down and when $Z>1$ the stock price jumps up. In the Merton jump diffusion model $Z=\exp(Y)$, in which $Y$ is a Gaussian random variable.
We will take $\mu=r+\lambda-\lambda \xi$, in which
$\xi=\int_{\R_+}x v(dx)<\infty$,
so that $(e^{-rt}S_t)_{t \geq 0}$ is a martingale, i.e., $\P$ is a risk neutral measure. The constant $r\geq 0$ is the interest
rate, and the constant $\sigma>0$ is the volatility.
 We assume the risk neutral pricing measure $\P$, and hence the parameters of the problem, are fixed as a result of a calibration to historical data.
The value function of the American put option pricing problem is
\begin{equation}\label{eq:opt-stop}
V(x,T):=\sup_{\tau \in \widetilde{\S}_{0,T}}\E\{e^{-\alpha
\tau}(K-S_{\tau})^+\},
\end{equation}
in which $\widetilde{\S}_{0,T}$ is the set of stopping times (of the
filtration generated by $W$ and $N$) that take values in $[0,T]$, and
$\E$ is the expectation under the probability measure $\P$, given
that $S_0=x$.

We will show that $V$ is the classical solution of a
free boundary equation and that it satisfies the \emph{smooth
fit principle}, i.e.,
$V$ is continuously differentiable with respect to
its first variable at the optimal stopping boundary. We argue these facts
by showing that $V$ is the fixed point of an
operator, which we will denote by $J$, that maps a given function  to the
value function of an optimal stopping problem for a
geometric Brownian motion. This operator acts as a regularizer: As soon as the given function $f$
has some certain regularity properties, we show that $Jf$ is
the unique classical solution of a corresponding free boundary equation
and that it satisfies the smooth fit principle. The proof of the main result concludes once we show that $V$ has these certain regularity
properties. In this last step we 
make use of a sequence (which is constructed by iterating $J$ starting with the
pay-off function of the put option) that converges to $V$
uniformly and exponentially fast. Incidentally, this sequence yields a
numerical procedure, whose accuracy versus speed characteristics
can be controlled. Each element of this sequence is an optimal stopping problems for
geometric Brownian motion and can be readily calculated using classical finite difference methods 
(see e.g. \cite{dewyne} for the implementation of these methods).
An alternative proof of the regularity of $V$ was given in
\cite{pham97}. This proof used a
combination of the results in \cite{friedmansde} and the theory of viscosity
solutions. In particular the proof of Proposition 3.1 in \cite{pham97} is carried out
 (details are not provided but hinted) using arguments similar to those used in the proof of Proposition 5.3 in \cite{pham1}. The latter proof
 uses the uniqueness results of \cite{Ishii} for viscosity solutions. 
  
 The infinite horizon American put option for jump diffusions were
analyzed in \cite{bayraktar-inf-horizon} using the iterative scheme we describe here. The
main technical difficulty in the current paper stems from the fact that each element in the approximating sequence solves a 
 parabolic rather an elliptic problem. 
In fact, in the infinite horizon case one can obtain a closed form representation for the value function, which is not possible in the finite horizon case. We make use of the results of \cite{friedmansde},
and Chapter 2 of \cite{karatzas-shreve-book2} (also see
 Chapter 7 of \cite{peskir-shiryaev-2006}) to study the properties
of the approximating sequence. For example, we show that the approximating sequence
is bounded with respect to the H\"{o}lder semi-norm
(see page 61 in \cite{friedman-parabolic} for a definition), which is used to argue that the limit
of the approximating sequence (which is a fixed point of $J$)
solves a corresponding free boundary equation.

Somewhat similar approximation techniques to the one we employ were used to solve
optimal stopping problems for \emph{diffusions}: see e.g. \cite{alvarez3} for perpetual optimal
stopping problems with non-smooth pay-off functions; and
\cite{carr}, \cite{touzi-karoui} for finite time horizon
American put option pricing problems for geometric Brownian motion. On the other hand, \cite{kyprianu} and \cite{sm} consider the smooth
fit principle for the infinite horizon American put option pricing problems for one-dimensional exponential
L\'{e}vy processes using the fluctuation theory. Also see \cite{bdk05} for the analysis of the smooth fit principle for a multi-dimensional infinite horizon optimal stopping problem.

The next two sections prepare the proof our main result, Theorems~\ref{thm:main}, in a sequence of lemmata and corollaries.
In
the next section, we introduce the functional operator $J$, that
maps a given function to the value function of an optimal stopping problem
for a geometric Brownian motion. We then analyze the properties of $J$. For example, $J$ preserves
convexity with respect to the first variable; the increase in the
H\"{o}lder semi-norm after the application of $J$ can be controlled; $J$ maps
certain class of functions to the classical solutions of
free boundary equations. In Section 3, we construct a
sequence of functions that converge to the smallest fixed point of the
operator $J$.
We show that the sequence is bounded in the H\"{o}lder norm,
and satisfies certain regularity properties using results of
Section 2. We eventually arrive at the fact that the smallest fixed point
of $J$ is equal to $V$. As a result the regularity properties of $V$ follow.
\section{A Functional Operator and Its Properties}
Let us define an operator $J$ through its action on a test
function $f:\R_+ \times \bar{\R}_+  \rightarrow \R_+$: The operator $J$ takes the function $f$ to the 
value function of the following optimal stopping problem
\begin{equation}\label{eq:defn-func-op}
J f (x,T)=\sup_{\tau \in
\S_{0,T}}\E\left\{\int_0^{\tau}e^{-(\alpha+\lambda)t}\lambda \cdot
Pf (S_t^0,T-t)dt+
e^{-(\alpha+\lambda)\tau}(K-S^0_{\tau})^+\right\},
\end{equation}
in which
\begin{equation}\label{eq:defn-Sf}
Pf (x,T-t)=\int_{\R_+}f(x z,T-t) \nu(dz), \quad x\geq 0.
\end{equation}
We will extend $T \rightarrow Jf(x,T)$ onto $[0,\infty]$ by
letting
\begin{equation}\label{eq:first-extension-eq}
Jf(x,\infty)=\lim_{T \rightarrow \infty}Jf(x,T).
\end{equation}
Here, $S^0=\{S^0_t;t \geq 0\}$ is the solution of
\begin{equation}
dS^0_t=\mu S^0_t dt+\sigma S_t^0 dW_t, \quad S^0_0=x,
\end{equation}
whose infinitesimal generator is given by
\begin{equation}\label{eq:inifinite}
\mathcal{A}:=\frac{1}{2}\sigma^2 x^2\frac{d^2}{dx^2}+\mu x
\frac{d}{dx}.
\end{equation}
In \eqref{eq:defn-func-op}, $\S_{[0,T]}$ denotes the set of stopping times of $S^0$ which take values in $[0,T]$.
Note that
\begin{equation}
S^0_t=x H_t,
\end{equation}
where
\begin{equation}
H_t=\exp\left\{\left(\mu-\frac{1}{2}\sigma^2\right)t+\sigma
W_t\right\}.
\end{equation}
The next remark characterizes the optimal stopping times of
(\ref{eq:defn-func-op}) using the Snell envelope theory.
\begin{remark}\label{rem:snell}
Let us denote
\begin{equation}
Y_t:=\int_0^{t}e^{-(\alpha+\lambda)s}\lambda \cdot Pf
(S_t^0,T-s)ds+ e^{-(\alpha+\lambda)t}(K-S^0_t)^+.
\end{equation}
Using the strong Markov property of $S^0$, we can determine the
Snell envelope of $Y$ as
\begin{equation}
\xi_t:=\sup_{\tau \in
\S_{t,T}}\mathbb{E}\left\{Y_{\tau}|\F_t\right\}=e^{-(\lambda+\alpha)t}Jf(S^0_t,T-t)+\int_0^{t}
e^{-(\alpha+\lambda)s} \lambda\;Pf(S^0_s,T-s)ds, \quad t \in
[0,T].
\end{equation}
Theorem D.12 in \cite{karatzas-shreve-book2} implies that
the stopping time
\begin{equation}\label{eq:opt-st-time}
\tau_x:=\inf\{t \in [0,T): \xi_t =Y_t\}\wedge T = \inf\{t \in
[0,T]: Jf(S_t^0,T-t)=(K-S^0_t)^+\},
\end{equation}
satisfies
\begin{equation}\label{eq:with-the-optimal-stopping-time}
Jf(x,T)=\E\left\{\int_0^{\tau_x}e^{-(\alpha+\lambda)t}\lambda
\cdot Pf (S_t^0,T-t)dt+
e^{-(\alpha+\lambda)\tau_x}(K-S^0_{\tau_x})^+\right\}.
\end{equation}
Moreover, the stopped process $(e^{-(r+\lambda)(t \wedge \tau_x)}Jf(S_{t \wedge \tau_x}^0,T-t\wedge \tau_x)+\int_0^{t \wedge \tau_x}e^{-(\alpha+\lambda)s}\lambda
\cdot Pf (S_s^0,T-s)ds)_{t \geq 0}$ is a martingale.
 The second infimum in (\ref{eq:opt-st-time}) is less
than $T$ because $Jf(S_T^0,0)=(K-S^0_T)^+$. 

When $f$ is bounded, it follows from the bounded convergence
theorem that (using the results of \cite{bayraktar-inf-horizon}
and arguments similar to the ones used in Corollary 7.3 in Chapter
2 of \cite{karatzas-shreve-book2})
\begin{equation}\label{eq:extension}
Jf(x,\infty)=\sup_{\tau \in \S_{0,\infty}
}\E\left\{\int_0^{\tau}e^{-(\alpha+\lambda)t}\lambda \cdot Pf
(S_t^0,\infty)dt+
e^{-(\alpha+\lambda)\tau}(K-S^0_{\tau})^+\right\}.
\end{equation}
\end{remark}

The next three lemmas on the properties of $J$ immediately follow from the definition in \eqref{eq:defn-func-op}. The first lemma states that $J$ preserves monotonicity.
\begin{lemma}\label{lem:mono}
Let $T \rightarrow f(x,T)$ be non-decreasing, and $x \rightarrow
f(x,T)$ be non-increasing. Then $T \rightarrow
Jf(x,T)$ is non-decreasing and $x \rightarrow Jf(x,T)$ is
non-increasing.
\end{lemma}
The operator $J$ preserves boundedness and
order.
\begin{lemma}\label{lem:pre-bnd}
Let $f:\R_+ \times \bar{\R}_+  \rightarrow \R_+$ be a bounded
function. Then $J f$ is also bounded. In fact,
\begin{equation}\label{eq:sup-norm-J-f}
0 \leq \|Jf\|_{\infty}\leq
K+\frac{\lambda}{\alpha+\lambda}\|f\|_{\infty}.
\end{equation}
\end{lemma}

\begin{lemma}\label{lem:preserve}
For any $f_1,f_2:\R_+ \times \bar{\R}_+  \rightarrow
\R_+$ that satisfy $f_1(x,T) \leq f_2(x,T)$, we have that $Jf_1 (x,T) \leq
Jf_2(x,T)$ for all $(x,T) \in \R_+ \times \bar{\R}_+ $ .
\end{lemma}

As we shall see next, the operator $J$ 
preserves convexity (with respect to the first variable).
\begin{lemma}\label{lem:convexity}
If $f:\R_+ \times
\bar{\R}_+  \rightarrow \R_+$ is a convex function in its first
variable, then so is $Jf: \R_+ \times \bar{\R}_+  \rightarrow
\R_+$.
\end{lemma}
\begin{proof}
Note that $Jf$ can be written as
\begin{equation}\label{eq:dmyJ}
J f (x,T)=\sup_{\tau \in
\S_{0,T}}\mathbb{E}\left\{\int_0^{\tau}e^{-(\alpha+\lambda)t}\lambda \cdot
Pf (x H_t,T-t)dt+
e^{-(\alpha+\lambda)\tau}(K-xH_{\tau})^+\right\}.
\end{equation}
 Since $f(\cdot,T-t)$ is convex, so is
$Pf(\cdot,T-t)$. As a result the integral with respect to time in
 \eqref{eq:dmyJ} is also convex in $x$. On the other hand, note that $(K-xH_{\tau})^+$ is also a convex function of $x$. Taking the expectation does not change the convexity with respect to $x$. Since the upper envelope
(supremum) of convex functions is convex, the result follows.
\end{proof}

\begin{remark}\label{rem:at-the-abs-bd}
Since $x=0$ is an absorbing boundary for the process $S^0$, for
any $f:\R_+ \times \bar{\R}_+  \rightarrow \R_+,
$
\begin{equation}
\begin{split}
J f(0,T)&=\sup_{t \in
\{0,T\}}\left\{\int_0^{t}e^{-(\alpha+\lambda)s} \lambda
f(0,T-s)ds+e^{-(\lambda+\alpha)t}K\right\}
\\&=\max\left\{K, \int_0^{T} e^{-(\alpha+\lambda)s}\lambda f(0,T-s)ds+
e^{-(\lambda+\alpha)T}K\right\}, \quad T \geq 0.
\end{split}
\end{equation}
If we further assume $f \leq K$, then $Jf (0,T) = K$,
$T \geq 0$. 
\end{remark}

\begin{lemma}\label{eq:J-f-lp}
Let us assume that $f:\R_+ \times \bar{\R}_+  \rightarrow \R_+$ is
convex in its first variable and $\|f\|_{\infty} \leq K$. Then
 $x \rightarrow Jf(x,t)$ satisfies
\begin{equation}\label{eq:Lips-x}
|Jf(x,T)-Jf(y,T)| \leq |x-y|, \quad (x,y) \in \R_+ \times
\bar{\R}_+ ,
\end{equation}
and all $T \geq 0$.
\end{lemma}
\begin{proof}
First note that a positive convex function that is bounded above has to be non-increasing. Therefore $f$ is non-increasing.
As a result of  Lemma \ref{lem:mono},
  $x \rightarrow Jf(x,t)$ is
non-increasing. This function is convex (by Lemma~\ref{lem:convexity}),
and it satisfies
\begin{equation}
Jf(x,T) \geq (K-x)^+ \quad Jf(0,T)=K.
\end{equation}
Consequently, the left and right derivatives of $Jf$ satisfy
\begin{equation}
-1 \leq D_{-}^x Jf(x,T) \leq D_+^xJf(x,T) \leq 0, \quad x>0,\; T
\geq 0.
\end{equation}
Now, the result follows since the derivatives are bounded by 1 (also see Theorem 24.7 (on page 237) in
\cite{rockafellar}).
\end{proof}
\begin{remark}\label{rem:pham-C}
Let $T_0 \in (0,\infty)$ and denote 
\begin{equation}\label{eq:amercan-dif}
F(x,T)= \sup_{\tau \in \S_{0 , T}}
\Eb\left\{e^{-(r+\lambda)\tau}(K-x H_{\tau})^+\right\}, \quad x
\in R_+,\; T \in [0,T_0].
\end{equation}
Then for $S \leq T \leq T_0$
\begin{equation}\label{eq:pham-s-ineq}
F(x,T)-F(x,S) \leq C \cdot |T-S|^{1/2},
\end{equation}
for all $x \in \R_+$ and for some $C$ that depends only on $T_0$.
See e.g.  equation (2.4) in \cite{pham97}.
\end{remark}
The next lemma, which is very crucial for our proof of the
smoothness of the American option price for jump diffusions, shows
that the increase in the H\"{o}lder semi-norm that the operator
$J$ causes can be controlled.

\begin{lemma}\label{lem:crucial}

Let us assume that for some $L \in (0,\infty)$
\begin{equation}\label{eq:holde}
|f(x,T)-f(x,S)| \leq L |T-S|^{1/2}, \quad (T,S) \in [S_0,T_0] \times [S_0,T_0],
\end{equation}
for all $x \in \R_+$, for $0 \leq S_0<T_0<\infty$. Then
\begin{equation}\label{eq:Jf-uni-holder}
|Jf(x,T)-Jf(x,S)| \leq (a \, L+C) \; |T-S|^{1/2}, \quad (T,S) \in
[S_0,T_0] \times [S_0,T_0],
\end{equation}
for some $a \in (0,1)$ whenever
\begin{equation}\label{eq:T}
|T-S| < \left(\frac{r}{r+\lambda} \frac{L}{\lambda K}\right)^2.
\end{equation}
Here, $C \in (0,\infty)$ is as in  Remark~\ref{rem:pham-C}.
\end{lemma}

\begin{proof}
 Without loss of generality
we will assume that $T>S$. Then we can write
\begin{equation}\label{eq:long-one}
\begin{split}
 &Jf(x,T)-Jf(x,S) \leq \sup_{\tau \in \S_{0,T}}
\Bigg[\Eb\left\{\int_0^{\tau }e^{-(r+\lambda)t}
\lambda\;Pf(xH_t,T-t)dt+e^{-(r+\lambda)\tau}
(K-xH_{\tau})^+\right\}
\\&-\Eb\left\{\int_0^{\tau \wedge S}e^{-(r+\lambda)t} \lambda\;Pf(xH_t,S-t)dt+e^{-(r+\lambda)(\tau \wedge S)} (K-xH_{\tau \wedge S
})^+\right\}\Bigg]
\\&=\sup_{\tau \in \S_{0,T}}
\Bigg[\Eb \bigg\{ \int_0^{\tau}e^{-(r+\lambda)t}
\lambda\,\left(Pf(xH_t,T-t)-Pf(xH_t,S-t)\right)dt
\\&+1_{\{S<\tau\}}\left[\int_{S}^{\tau}e^{-(r+\lambda)t}
\lambda\, Pf(xH_t,S-t)dt+
\left(e^{-(r+\lambda)\tau}(K-xH_{\tau})^+-e^{-(r+\lambda)S}(K-xH_{S})^+\right)\right]\bigg\}\Bigg]
\\ &\leq \frac{\lambda}{r+\lambda} L \; (T-S)^{1/2}+\frac{\lambda}{r+\lambda} \,
K\,\left(e^{-(r+\lambda)S}-e^{-(r+\lambda)T}\right)
\\& \qquad \qquad+\sup_{\tau \in \S_{S,T}}\Eb\left\{\left(e^{-(r+\lambda)\tau}(K-xH_{\tau})^+\right\}-\Eb\left\{e^{-(r+\lambda)S}(K-xH_{S})^+\right)\right\}
\\&\leq \frac{\lambda}{r+\lambda} L \; (T-S)^{1/2}+ \lambda \, K
\;(T-S)+e^{-(r+\lambda)S}\left(F(H_S,T-S)-F(H_S,0)\right),
\\& \leq  \left(\frac{\lambda}{r+\lambda} L + C\right) \; (T-S)^{1/2}+ \lambda \, K
\;(T-S)
\end{split}
\end{equation}
in which $F$ is given by (\ref{eq:amercan-dif}). To
derive the second inequality in \eqref{eq:long-one}, we use the fact that
\begin{equation}
\begin{split}
 |Pf(xH_t,T-t)-Pf(xH_t, S-t)| &\leq \int_{\R_+} \nu(dz) \left|f(x
z H_t,T-t)-f(x zH_t, S-t)\right| \leq  L \; |T-S|^{1/2},
\end{split}
\end{equation}
which follows from the assumption in (\ref{eq:holde}), and that
\begin{equation}
\Eb \left\{1_{\{S<\tau\}}\int_0^{\tau \wedge S}e^{-(r+\lambda)t}
\lambda\;P f(xH_t,S-t)dt\right\}\leq \lambda K \Eb \left\{\int_{S}^{T}
e^{-(r+\lambda)t}dt \right\}\leq
\frac{\lambda\,K}{\lambda+K}
\left(e^{-(r+\lambda)S}-e^{-(r+\lambda)T}\right).
\end{equation}
To derive the third inequality in \eqref{eq:long-one}, we use
\begin{equation}
e^{-(r+\lambda)S}-e^{-(r+\lambda)T} \leq
e^{-(r+\lambda)S}(r+\lambda)(T-S) \leq (r+\lambda)(T-S).
\end{equation}
The last inequality in \eqref{eq:long-one} follows from (\ref{eq:pham-s-ineq}).
Equation
(\ref{eq:Jf-uni-holder}) follows from (\ref{eq:long-one}) whenever $T$ and $S$ satisfy (\ref{eq:T}).
\end{proof}
\noindent Let us define the continuation
region and its sections by
\begin{equation}\label{eq:cont-reg}
\C^{Jf}:=\{(T,x) \in (0,\infty)^2: Jf(x,T)>(K-x)^+\}, \;
\text{and} \; \C^{Jf}_T:=\{x \in (0,\infty): Jf(T,x)>(K-x)^+\},\;
\end{equation}
$T>0$, respectively.

\begin{lemma}\label{lem:cont-reg}
Suppose that $f:\R_+ \times \bar{\R}_+ \rightarrow \bar{\R}_+$ is such
that $x \rightarrow f(x,T)$ is a positive convex function, $T
\rightarrow f(x,T)$ is non-decreasing, and $\|f\|_{\infty}\leq K$.
Then for every $T>0$ there exists $c^{Jf}(T) \in (0,K)$ such that
$\C^{Jf}_T=(c^{Jf}(T),\infty)$. Moreover, $T \rightarrow c^{Jf}(T)$ is
non-increasing.
\end{lemma}

\begin{proof}
Let us first show that if $x \geq K$, then $x \in \C^{Jf}_{T}$
for all $T \geq 0$. Let
$\tau_{\eps}:=\inf\{0 \leq t \leq T: S_t^0 \leq K-\eps\}$.
Since
$
\P\{0<\tau_{\eps}<T\}>0 \quad \text{for} \; x \geq K$,
for all $T> 0$, we have that
\begin{equation}\label{eq:eps-geq-K}
\E \left\{\int_0^{\tau_{\eps}}e^{-(r+\lambda)t}
\lambda\;Pf(S_t^0,T-t)dt+e^{-(r+\lambda)\tau_{\eps}}
(K-S^0_{\tau_{\eps}})^+\right\} >0,
\end{equation}
which implies that $x \in \C^{Jf}_T$. 
On the other hand, it is clear that
\begin{equation}\label{eq:trivia}
(K-x)^+ \leq Jf(x,T) \leq J f(x,\infty), \quad (x,T) \in \R_+
\times \bar{\R}_+ .
\end{equation}
 Thanks to
in Lemma 2.6 of \cite{bayraktar-inf-horizon}, there exist $l^f \in
(0,K)$ such that
\begin{equation}\label{eq:hat-j-toches}
J f(x,\infty)=(K-x)^+, \; x \in [0,l^f]; \quad J
f(x,\infty)>(K-x)^+, \; x \in (l^f,\infty).
\end{equation}
Since $x \rightarrow J f(x,\infty)$ and $x \rightarrow Jf(x,T)$, $T \geq
0$, are convex functions (from Lemma 2.2 in
\cite{bayraktar-inf-horizon} and Lemma~\ref{lem:convexity}
respectively), (\ref{eq:eps-geq-K}), (\ref{eq:trivia}) and
(\ref{eq:hat-j-toches}) imply that there exists a point $c^{Jf}(T) \in
(l^f,K)$ such that
\begin{equation}\label{eq:Jf-l-f}
J f(x)=(K-x)^+, \; x \in [0,c^{Jf}(T)]; \quad J f(x,T)>(K-x)^+, \; x
\in (c^{Jf}(T),\infty),
\end{equation}
for $T>0$. This proves the first statement of the Lemma.
The fact that $T \rightarrow c(T)$ is non-increasing follows from
the fact that $T \rightarrow Jf(x,T)$ is non-decreasing.
\end{proof}

In the following lemma we will argue that if $f$ has
certain regularity properties, then $Jf$ is the
classical solution of a parabolic free boundary equation.

\begin{lemma}\label{eq:lem-KS}
Let us assume that $f:\R_+ \times \bar{\R}_+ \rightarrow \R_+$ is
convex in its first variable, $\|f\|_{\infty} \leq K$ and $T
\rightarrow f(x,T)$ is non-increasing. Moreover, we will assume that $f$ satisfies
\begin{equation}
|f(x,T)-f(x,S)| \leq A\, |T-S|^{1/2} \quad \text{whenever}\;\;\; |T-S|<B,
\end{equation}
for all $x \in \R_+$, where $A, B$ are strictly positive constants that do not depend on $x$. Then the function
$Jf:\R_+ \times \R_+  \rightarrow \R_+$ is the unique bounded solution (in the classical sense) of 
\begin{align}
&\A u(x,T)- (\alpha+\lambda)\cdot u (x,T)+\lambda \cdot (Pf)(x,T)
- \frac{\partial}{\partial T} u(x,T)=0 \quad x>c^{Jf}(T), \label{eq:cont-regp}
\\ & u(x,T) =(K-x) \quad x \leq c^{Jf}(T), \label{eq:stop-reg}
\end{align}
in which $\A$ is as in (\ref{eq:inifinite}) and $c^{Jf}$ is as in
Lemma~\ref{lem:cont-reg}. Moreover,
\begin{equation}\label{eq:thnkstoHao}
\A Jf(x,T)- (\alpha+\lambda)\cdot Jf (x,T)+\lambda \cdot (Pf)(x,T)
- \frac{\partial}{\partial T} Jf(x,T) \leq 0, \quad x<c^{Jf}(T).
\end{equation}
\end{lemma}

\begin{proof}
The proof is motivated by Theorem 2.7.7
of \cite{karatzas-shreve-book2}. 
Equation (\ref{eq:stop-reg}) is clearly satisfied by $Jf$. In what follows, we will first show that $Jf$ satisfies (\ref{eq:cont-regp}). Let us take a point in $(t,T) \in \mathcal{C}^{Jf}$ and
consider a bounded rectangle $R=(t_1,t_2) \times (x_1,x_2)$ containing this point. We will let
\begin{equation}
t_2-t_1<B \wedge \left(\frac{r A}{(r+\lambda) \lambda K}\right)^2.
\end{equation}
Let
$\partial_0 R$ be the parabolic boundary of $R$
and consider the parabolic partial differential equation 
\begin{equation}\label{eq:parabol}
\begin{split}
&\A u(x,T)- (\alpha+\lambda)\cdot u (x,T)+\lambda \cdot (Pf)(x,T)
- \frac{\partial}{\partial T} u(x,T)=0 \quad \text{in}\, R,
\\ &u(x,T)=Jf(x,T) \quad \text{on $\partial_0 R$}.
\end{split}
\end{equation}
As a result of Lemmas \ref{eq:J-f-lp} and \ref{lem:crucial}, $Jf$ satisfies the uniform Lipschitz and
H\"{o}lder continuity conditions, which implies that $Jf$ is continuous. On the other hand, for any $(T,x) \in R$
\begin{equation}
\begin{split}
|Pf(x,T)-Pf(y,S)| &\leq |Pf(x,T)-Pf(x,S)|+|Pf(x,S)-Pf(y,S)|
\\& \leq \int_{\R_+} \nu(dz) \left(|f(x z,T)-f(xz, S)|+ |f(xz,S)-f(yz,S)|\right)
\\&\leq A \; |T-S|^{1/2}+ \xi\; |x-y|,
\end{split}
\end{equation}
Now, Theorem 5.2 in \cite{friedmansde} implies that (\ref{eq:parabol}) has a unique classical solution.
We will show that this unique solution coincides with $Jf$ using optional sampling theorem. Let us introduce the stopping time
\begin{equation}
\tau:=\inf\{ \theta \in [0,t_0-t_1): (t_0-\theta, x_0 H_\theta \in \partial_0 R\} \wedge (t_0-t_1),
\end{equation}
which is the first time $S^0$ hits the parabolic boundary when $S^0$ starts from $(x_0,t_0)$. Let us also define the process
$N_{\theta}:=e^{-(r+\lambda) \theta} u(x_0 H_\theta,t_0-\theta)+\int_0^{\theta}e^{-(\alpha+\lambda)t}\lambda
\cdot Pf (S_t^0,t_0-t)dt$, $\theta \in [0, t_0-t_1]$. From the classical It\^{o}'s formula it follows that the stopped process $N_{\theta \wedge \tau}$ is a bounded martingale. As a result
\begin{equation}\label{eq:uxoto}
u(x_0,t_0)=N_0=\E \left\{N_{\tau}\right\}=\mathbb{E}\left\{e^{-(r+\lambda) \tau} Jf(x H_\tau, t_0-\tau))+\int_0^{\tau}e^{-(\alpha+\lambda)t}\lambda
\cdot Pf (S_t^0,t_0-t)dt\right\}.
\end{equation}
Clearly $\tau \leq \tau_x$. Since the stopped process $(e^{-(r+\lambda)(t \wedge \tau_x)}Jf(S_{t \wedge \tau_x}^0,t_0-t\wedge \tau_x)+\int_0^{t \wedge \tau_x}e^{-(\alpha+\lambda)t}\lambda
\cdot Pf (S_s^0,t_0-s)ds)_{t \geq 0}$ is a bounded martingale, another application of the optional sampling theorem yields
\begin{equation}\label{eq:Jfxoto}
\mathbb{E}\left\{e^{-(r+\lambda) \tau}Jf(x_0 H_\tau, t_0-\tau)+\int_0^{\tau}e^{-(\alpha+\lambda)t}\lambda
\cdot Pf (S_t^0,t_0-t)dt\right\}=Jf(x_0,t_0).
\end{equation}
Combining (\ref{eq:uxoto}) and (\ref{eq:Jfxoto}), we see that (\ref{eq:cont-regp}) is satisfied in the classical sense since the choice of $(x_0,t_0) \in \C^{Jf}$ is arbitrary.

We still need to show uniqueness among bounded functions. Fix $x>c^{Jf}(T)$. Let $u$ be a bounded function satisfying 
(\ref{eq:cont-regp}) and (\ref{eq:stop-reg}). Let us define $M_t:=e^{-(r+\lambda)t} u(xH_t,T-t)+\int_0^{t}e^{-(\alpha+\lambda)t}\lambda
\cdot Pf (S_s^0,T-s)ds$. Using the classical It\^{o} formula it can be seen that $M_{t \wedge \tau_x}$ is a bounded martingale. Since $\tau_x$ is optimal (see \eqref{eq:with-the-optimal-stopping-time}), by the optional sampling theorem, we have 
\begin{equation}
\begin{split}
u(x,T)&=M_0=\E\{M_{\tau_x}\}=\mathbb{E}\left\{e^{-(r+\lambda) \tau}u(x H_{\tau_x}, T-\tau_x)+\int_0^{\tau_x}e^{-(\alpha+\lambda)t}\lambda
\cdot Pf (S_s^0,T-s)ds\right\}
\\&=\mathbb{E}\left\{e^{-(r+\lambda) \tau}(K-xH_{\tau_x})^++\int_0^{\tau_x}e^{-(\alpha+\lambda)t}\lambda
\cdot Pf (S_s^0,T-s)ds\right\}=Jf(x,T).
\end{split}
\end{equation} 
Next, we will prove \eqref{eq:thnkstoHao}. To this end, let $x<c^{Jf}(t)$. Let $U$ a closed interval centered at $x$ such that $U \subset (0,c^{Jf}(T))$. Let $\tau_U=\{t \geq 0: x H_{t} \notin U\}$. Since $(e^{-(r+\lambda)t}Jf(S_{t }^0,T-t)+\int_0^{t }e^{-(\alpha+\lambda)s}\lambda
\cdot Pf (S_s^0,T-s)ds)_{t \geq 0}$ is a supermartingale we can write
\begin{equation}
\mathbb{E} \left[e^{-(r+\lambda)(\tau_U \wedge t)}Jf(x H_{\tau_U \wedge t}, T-\tau_U \wedge t))+\int_0^{\tau_U \wedge t}e^{-(r+\lambda)u} \lambda Pf(xH_u,T-u)du\right] \leq Jf(x,T),
\end{equation}
for all $t \geq 0$.
Since $Jf(x,t)=K-x$ when $(T,x) \in \R_+^2-C^{Jf}$, we can apply It\^{o}'s formula to obtain that
\begin{equation}
\lim_{t \rightarrow 0}\mathbb{E}\left[\frac{1}{t}\int_0^{\tau_U \wedge t}e^{-(r+\lambda)u}\left(\left(\A-(r+\lambda)\cdot-\frac{\partial}{\partial T} \right) Jf(xH_u,t-u)+\lambda Pf(xH_u,T-u) \right)du\right] \leq 0.
\end{equation}
 Now, \eqref{eq:thnkstoHao} follows thanks to dominated convergence theorem, which allows us to exchange the limit and the expectation. We can apply the dominated convergence theorem thanks to the fact that $U$ is a compact domain. 
\end{proof}
\begin{lemma}\label{lem:smooth-fit-cond}
For a given $T>0$, let $x \rightarrow f(x,T)$ be a convex and
non-increasing function. Then the convex function $x \rightarrow
Jf(x,T)$ is of class $C^1$ at $x=c(T)$, i.e.,
\begin{equation}
\frac{\partial}{\partial x} Jf(x,T)\bigg|_{x=c(T)}=-1.
\end{equation}

\end{lemma}

\begin{proof}
The proof is similar to the proof of Lemma 7.8 on page 74 of
\cite{karatzas-shreve-book2}, but we will provide it here for the
sake of completeness. If we let $x=c(T)$, then
\begin{equation}\label{eq:half-pf}
\begin{split}
Jf(x+\eps,T)&=\mathbb{E}\left\{\int_0^{\tau_{x+\eps}}e^{-(\alpha+\lambda)t}\lambda
\cdot Pf ((x+\eps)H_t,T-t)dt+
e^{-(\alpha+\lambda)\tau_{x+\eps}}(K-(x+\eps)H_{\tau_{x+\eps}})^+\right\}
\\ &=\mathbb{E}\left\{\int_0^{\tau_{x+\eps}}e^{-(\alpha+\lambda)t}\lambda \cdot
Pf (xH_t,T-t)dt+ e^{-(r+\lambda \tau_{x+\eps})} (K-x H_{\tau_{x+\eps}})^+\right\}
\\&+\mathbb{E}\left\{\int_0^{\tau_{x+\eps}}e^{-(\alpha+\lambda)t}\lambda \cdot
\left[Pf ((x+\eps)H_t,T-t)-Pf (x H_t,T-t)\right]dt\right\}
\\& +\mathbb{E}\left\{e^{-(r+\lambda)\tau_{x+\eps}}\left[(K-(x+\eps)H_{\tau_{x+\eps}})^+-(K- x H_{\tau_{x+\eps}})^+\right]\right\}
\\ &\leq Jf(x,T)+\mathbb{E}\left\{1_{\{\tau_{x+\eps}<T\}}e^{-(r+\lambda)\tau_{x+\eps}}\left[(K-(x+\eps)H_{\tau_{x+\eps}})-(K- x H_{\tau_{x+\eps}})\right]\right\}
\\&+\mathbb{E}\left\{1_{\{\tau_{x+\eps}=T\}}e^{-(r+\lambda)\tau_{x+\eps}}\left[(K-(x+\eps)H_{\tau_{x+\eps}})^+-(K- x H_{\tau_{x+\eps}})^+\right]\right\}
\\& \leq Jf(x,T)-\eps \E\left\{1_{\{\tau_{x+\eps}<T\}}e^{-(r+\lambda)\tau_{x+\eps}} H_{\tau_{x+\eps}}\right\}
\\&=Jf(x,T)-\eps \E\left\{e^{-(r+\lambda)\tau_{x+\eps}} H_{\tau_{x+\eps}}\right\}+\eps \E\left\{1_{\{\tau_{x+\eps}=T\}}e^{-(r+\lambda)T} H_{T}\right\}.
\end{split}
\end{equation}
The first inequality follows since $\tau_{x+\eps}$ is not optimal
when $S^0$ starts at $x$ and $x \rightarrow Pf(x,T)$ is a
decreasing function for any $T \geq 0$.  From (\ref{eq:half-pf}) it follows that
\begin{equation}
D_{+}^x Jf(x+,T) \leq -1,
\end{equation}
since $e^{-(\alpha+\lambda)t}H_t$ is a uniformly integrable
martingale and  $\tau_{x+\eps} \downarrow 0$. Convexity of
$Jf(t,x)$ (Lemma~\ref{lem:convexity}) implies that
\begin{equation}
-1=  D_{-}^{x}Jf(x-,t) \leq
 D_+^x Jf (x+,t) \leq -1,
\end{equation}
which yields the desired result.
\end{proof}

\section{A Sequence of Functions Approximating $V$}
\noindent Let us define a sequence of functions by the following iteration:
\begin{equation}\label{eq:defn-v-n}
v_0(x,T)=(K-x)^+, \quad v_{n+1}(x,T)=J v_{n}(x,T),\; n \geq 0,
\quad \text{for all $(x,T) \in \R_+ \times \R_+ $.}
\end{equation}
We extend these functions onto $\R_+ \times \bar{\R}_+$ by letting
\begin{equation}
v_{n}(x,\infty)=\lim_{T \rightarrow \infty}v_n(x,T).
\end{equation}
This sequence of functions is a bounded sequence as the next lemma
shows.
\begin{cor}\label{cor:bounded-seq}
For all $n
\geq 0$,
\begin{equation}\label{eq:sup-norm-v-n}
(K-x)^+ \leq v_{n}(x,T) \leq
\left(1+\frac{\lambda}{\alpha}\right)K, \quad (x,T) \in\R_+ \times
\bar{\R}_+ .
\end{equation}
\end{cor}
\begin{proof}
The first inequality follows since it may not be optimal to stop
immediately. Let us prove the second inequality using an induction
argument: Observe that $v_0(x,T)=(K-x)^+$, $(x,T) \in \R_+ \times
\bar{\R}_+ $, satisfies (\ref{eq:sup-norm-v-n}). Let us assume that
(\ref{eq:sup-norm-v-n}) holds for $n$ and show that it
holds for $n+1$. Using (\ref{eq:sup-norm-J-f}), we get that
\begin{equation}
\|v_{n+1}\|_{\infty}=\|Jv_{n}\|_{\infty}\leq
K+\frac{\lambda}{\alpha+\lambda}\left(1+\frac{\lambda}{\alpha}\right)K=\left(1+\frac{\lambda}{\alpha}\right).
\end{equation}
\end{proof}

As a corollary of Lemmas~\ref{lem:preserve} and \ref{lem:convexity}
 we can state the following corollary, whose
proof can be carried out by induction.
\begin{cor}\label{cor:increas}
The sequence $(v_{n}(x,T))_{n\geq 0}$ is
increasing for all $(x,T) \in \R_+ \times \bar{\R}_+ $. For each $n$, the
function $x \rightarrow v_n(x,T)$, $x \geq 0$, is convex for all
$T \in \bar{\R}_+$.
\end{cor}

\begin{remark}\label{rem:v-infty}
Let us define,
\begin{equation}\label{eq:defn-v-infty}
v_{\infty}(x,T):=\sup_{n \geq 0}v_n(x,T), \quad (x,T) \in\R_+
\times \bar{\R}_+ .
\end{equation}
This function is well defined as a result of
(\ref{eq:sup-norm-v-n}) and Corollary~\ref{cor:increas}. In fact,
it is convex, because it is the upper envelope of convex functions, and it is bounded by the right-hand-side of
(\ref{eq:sup-norm-v-n}).
\end{remark}

\begin{cor}\label{cor:decreasing}
For each $n \geq 0$ and $t \in \R_+$, $x \rightarrow v_n(x,T)$, is
a decreasing function on $[0,\infty)$. Moreover, $T \rightarrow
v_n(x,T)$ is non-decreasing. The same statements hold for $x \rightarrow
v_{\infty}(x,T)$,  and $T \to v_{\infty}(x,T)$, respectively.
\end{cor}

\begin{proof}
The behaviour with respect to the first variable is a result of
Corollary~\ref{cor:increas} and Remark~\ref{rem:v-infty} since any
positive convex function that is bounded from above is decreasing. 
For each $n$, the fact that $T \to v_n(x,T)$ is non-decreasing is a corollary of Lemma~\ref{lem:mono}. On the other hand, for any $T \geq S \geq 0$, we have that
$
v_{\infty}(x,T) =\sup_{n}v_{n}(x,T) \geq \sup_n v_{n} (x,S)=v_{\infty}(x,S)$.
\end{proof}
Next, we will sharpen the upper bound in
Corollary~\ref{cor:bounded-seq}. This improvement has some implications
for the continuity of $x \rightarrow v_{n}(x,T)$, $n \geq 1$,
and $x \rightarrow v_{\infty}(x,T)$ at $x=0$.
\begin{remark}\label{rem:sharp-upper-bound}
 The upper bound in
(\ref{cor:bounded-seq}) can be sharpened using
Corollary~\ref{cor:decreasing} and Remark~\ref{rem:at-the-abs-bd}.
Indeed, we have
\begin{equation}
(K-x)^+ \leq v_{n}(x,T) < K, \quad \text{for each $n$}, \quad
\text{and} \quad (K-x)^+ \leq v_{\infty}(x,T) <K,  \quad (x,T) \in
(0,\infty)^2.
\end{equation}
It follows from this observation that for every $T \in \bar{\R}_+$,
$x \rightarrow v_{n}(x,T)$, for every $n$, and $x \rightarrow
v_{\infty}(x,T)$, are continuous at $x=0$ since $v_{n}(0,T)=v_{\infty}(0,T)=K$ and these functions are convex. (Note that convexity already guarantees continuity for $x>0$.)
\end{remark}
\begin{lemma}\label{lem:fixed-po}
The function $v_{\infty}$ is the smallest fixed point
of the operator $J$.
\end{lemma}
\begin{proof}
\begin{equation}
\begin{split}
v_{\infty}(x,T-t)&=\sup_{n\geq 1}v_{n}(x,T-t)
\\&=\sup_{n \geq
1}\sup_{\tau \in \S_{0,T}
}\E\left\{\int_0^{\tau}e^{-(\alpha+\lambda)t}\lambda \cdot P v_n
(S_t^0,T-t)dt+ e^{-(\alpha+\lambda)\tau}(K-S^0_{\tau})^+\right\}
\\&=\sup_{\tau
\in \S_{0,T} }\sup_{n \geq
1}\E\left\{\int_0^{\tau}e^{-(\alpha+\lambda)t}\lambda \cdot P v_n
(S_t^0,T-t)dt+ e^{-(\alpha+\lambda)\tau} (K-S^0_{\tau})^+\right\}
\\&=\sup_{\tau
\in \S_{0,T} }\E\left\{\int_0^{\tau}e^{-(\alpha+\lambda)t}\lambda
\cdot P (\sup_{n \geq 1} v_n) (S_t^0,T-t)dt+
e^{-(\alpha+\lambda)\tau} (K-S^0_{\tau})^+\right\}
\\&=J v_{\infty}(x,T-t),
\end{split}
\end{equation}
in which the fourth equality follows by applying the monotone convergence
theorem three times. Let $w:\R_+ \times \bar{\R}_+  \rightarrow \R_+$ be
another fixed point of the operator $J$. We will argue by induction that $w \geq v_{\infty}$. For $(x,t) \in \R_+
\times \bar{\R}_+ $, $w(x,T-t)=Jw(x,T-t)$, which implies that
$w(x,T-t)=Jw(x,T-t) \geq (K-x)^+=v_0(\cdot)$. If we assume that
$w(x,T-t) \geq v_{n}(x,T-t)$, then $w(x,T-t) = Jw(x,T-t) \geq J
v_n(x,T-t)=v_{n+1}(x,T-t)$. Consequently $w(x,T-t) \geq v_n(x,T-t)$
for all $n \geq 0$. As a result $w(x,T-t) \geq \sup_{n\geq 0}
v_{n}(x,T-t) = v_{\infty}(x,T-t)$.
\end{proof}

\begin{lemma}\label{lem:fix}
The sequence $\{v_n(\cdot, \cdot)\}_{n \geq 0}$ converges
uniformly to $v_{\infty}$. In fact, the rate of
convergence is exponential:
\begin{equation}\label{eq:iniform}
v_{n}(x,T) \leq v_{\infty}(x,T) \leq
v_{n}(x,T)+\left(\frac{\lambda}{\lambda+\alpha}\right)^n K, \quad
(x,T) \in\R_+ \times \bar{\R}_+ .
\end{equation}
\end{lemma}

\begin{proof}
The first inequality follows from the definition of
$v_{\infty}$. The second inequality can be proved by
induction. The inequality holds when we set $n=0$ by
Remark~\ref{rem:sharp-upper-bound}. Assume that the inequality
holds for $n>0$. Then
\begin{equation}\label{eq:con-rate}
\begin{split}
&v_{\infty}(x,T)=\sup_{\tau \in \S_{0,T}
}\E\left\{\int_0^{\tau}e^{-(\alpha+\lambda)t}\lambda \cdot P
v_{\infty} (S_t^0,T-t)dt+ e^{-(\alpha+\lambda)\tau}
(K-S^0_{\tau})^+\right\}
\\& \leq \sup_{\tau \in \S_{0,T}
}\E\bigg\{\int_0^{\tau}e^{-(\alpha+\lambda)t}\lambda \cdot P v_n
(S_t^0,T-t)dt+ e^{-(\alpha+\lambda)\tau} (K-S^0_\tau)^+\bigg\}
+\int_0^{\infty}dt\,
e^{-(\lambda+\alpha)t}\lambda\left(\frac{\lambda}{\lambda+\alpha}\right)^n
K 
\\&=v_{n+1}(x,T)+\left(\frac{\lambda}{\lambda+\alpha}\right)^{n+1}
K.
\end{split}
\end{equation}
\end{proof}
\begin{remark}
Note that, for a fixed $T_0>0$,
\begin{equation}\label{eq:conv-rate}
v_{n}(x,T) \leq v_{\infty}(x,T) \leq
v_{n}(x,T)+\left(1-e^{-(r+\lambda)T_0}\right)^n\left(\frac{\lambda}{\lambda+\alpha}\right)^n
K, \quad x \in \R_+,\; T \in (0,T_0).
\end{equation}
This can be derived using an induction argument similar to the one used in the proof of Lemma~\ref{lem:fix}. We simply replace (\ref{eq:con-rate}) by
\begin{equation}
\begin{split}
&v_{\infty}(x,T)\leq \sup_{\tau \in \S_{0,T}
}\E\bigg\{\int_0^{\tau}e^{-(\alpha+\lambda)t}\lambda \cdot P v_n
(S_t^0,T-t)dt+ e^{-(\alpha+\lambda)\tau} (K-S^0_\tau)^+\bigg\}
\\&
+\int_0^{T_0}dt\,
e^{-(\lambda+\alpha)t}\left(1-e^{-(r+\lambda)T_0}\right)^n\lambda\left(\frac{\lambda}{\lambda+\alpha}\right)^n
K =v_{n+1}(x,T)+K \left(1-e^{-(r+\lambda)T_0}\right)^{n+1}\left(\frac{\lambda}{\lambda+\alpha}\right)^{n+1}
\end{split}
\end{equation}
Observe that one can replace $K$ in (\ref{eq:conv-rate}) by $\|v_{\infty}-v_0\|_{\infty}$.
Note that the convergence rate in (\ref{eq:conv-rate}) is fast. This will lead to a
numerical scheme, whose error
versus accuracy characteristics can be controlled, for pricing American options. 
\end{remark}

\begin{remark}\label{rem:v-1}
Let $T_0 \in (0,\infty)$. It can be shown using similar
arguments to the ones used in the proof of Lemma~\ref{lem:crucial}
that 
\begin{equation}
|v_{1}(x,T)-v_1(x,S)| \leq \lambda\, K\; |T-S|+ C \, |T-S|^{1/2}, \quad T,S \in (0,T_0],
\end{equation}
for all $x \in \R_+$, in which  $C \in (0,\infty)$ is as in
Remark~\ref{rem:pham-C}. In fact
\begin{equation}
|v_{1}(x,T)-v_1(x,S)| \leq L \, |T-S|^{1/2},
\end{equation}
for all $x \in\R_+$ and for some $L$ that depends only on $T_0$.
\end{remark}

The next lemma shows that the functions $v_{n}$, $n \geq 0$, and $v_{\infty}$ are locally H\"{o}lder continuous with respect to the time variable.

\begin{lemma}\label{eq:v-infty-n-hold}
Let $T_0 \in (0,\infty)$ and $L \in (0,\infty)$ be as in
Remark~\ref{rem:v-1} and $C \in (0,\infty)$ be as in
Remark~\ref{rem:pham-C}. Then for $T,S \in (0,T_0)$, we have that
\begin{equation}\label{eq:ind}
|v_{n}(x,T)-v_n(x,S)| \leq \left(L+\frac{C}{1-a}
\right)|T-S|^{1/2} \quad \text{whenever} \;\;\; |T-S| \leq
\left(\frac{r}{r+\lambda}\;\frac{L}{\lambda \, K}\right)^2,
\end{equation}
for all $x \in \R_+$ and for all $n \geq 1$. Here, $a \in (0,1)$
is as in Lemma~\ref{lem:crucial}. Moreover,
\begin{equation}\label{eq:hold-v-infty}
|v_{\infty}(x,T)-v_{\infty}(x,S)| \leq  \left(L+\frac{C}{1-a}
\right) |T-S|^{1/2} \quad \text{whenever} \;\;\; |T-S| \leq
\left(\frac{r}{r+\lambda}\;\frac{L}{\lambda \, K}\right)^2,
\end{equation}
for all $x \in \R_+$.
\end{lemma}

\begin{proof}
The proof of (\ref{eq:ind}) will be carried out using an induction
argument. Observe from Remark~\ref{rem:v-1} that (\ref{eq:ind})
holds for $n=1$. Let us assume that (\ref{eq:ind}) holds for $n$
and show that it holds for $n+1$. Using Lemma~\ref{lem:crucial}, we
have that
\begin{equation}\label{eq:indunction}
|v_{n+1}(x,T)-v_{n+1}(x,S)| \leq \left(a \left(L+ \frac
{C}{1-a}\right)+C\right)|T-S|^{1/2},
\end{equation}
for  $|T-S| \leq
\left(\frac{r}{r+\lambda}\; \frac{L+C/(1-a)}{\lambda\,K}\right)^2$.
It is clear that the right-hand-side of (\ref{eq:indunction}) is
less than that of (\ref{eq:ind}), and
\begin{equation}
\frac{r}{r+\lambda}\; \frac{L+C/(1-a)}{\lambda\,K} \geq
\frac{r}{r+\lambda}\;\frac{L}{\lambda \, K},
\end{equation}
from which the first statement of the lemma follows.
Now let us prove (\ref{eq:hold-v-infty}). To this end observe that
\begin{equation}\label{eq:hold-v-inf-with-n}
\begin{split}
|v_{\infty}(x,T)-v_{\infty}(x,S)| &\leq
|v_{\infty}(x,T)-v_{n}(x,T)|+|v_{n}(x,T)-v_{n}(x,S)|+|v_{\infty}(x,S)-v_{n}(x,S)|
\\& \leq 2 \left(\frac{\lambda}{\lambda+r}\right)^n \, K+ \left(L+\frac{C}{1-a}\right)|T-S|^{1/2},
\end{split}
\end{equation}
for any $n>1$, which follows from (\ref{eq:ind}) and
Lemma~\ref{lem:fix}. The result follows since $n$ on the right-hand-side of (\ref{eq:hold-v-inf-with-n}) is arbitrary.
\end{proof}

\begin{lemma}
For $n \geq 0$, 
$|v_{n}(x,T)-v_n(y,T)| \leq |x-y|$, and
$|v_{\infty}(x,T)-v_{\infty}(y,T)| \leq |x-y|$, $(x,y) \in \R_+
\times \bar{\R}_+$, 
for all $T \geq 0$.
\end{lemma}
\begin{proof}
It follows from Remark~\ref{rem:sharp-upper-bound} that
$\|v_{n}\|_{\infty} \leq K$, for all $n \geq 0$, and
$\|v_{\infty}\|_{\infty} \leq K$. Moreover, for each $n \geq 0$,
$v_{n}(\cdot,T)$ is convex (for all $T \in \bar{\R}_+$) as a
result of Corollary~\ref{cor:increas}. On the other hand, it was
pointed in Remark~\ref{rem:v-infty} that $v_{\infty}(\cdot,T)$ is
convex for all $T \in \R_+$. Since
\begin{equation}\label{eq:jvs}
v_{n+1}(x,T)=Jv_{n}(x,T) \quad \text{and} \quad
v_{\infty}(x,T)=Jv_{\infty}(x,T),
\end{equation}
the statement of the lemma follows from Lemma~\ref{eq:J-f-lp}. 
\end{proof}

\begin{lemma}\label{lem:before-main}
For all $T \geq 0$ and $ n\geq 0$,
$\C^{v_{n+1}}_T=(c^{v_{n+1}}(T),\infty)$ for some $c^{v_{n+1}}(T) \in
(0,K)$ and $\C^{v_{\infty}}_T=(c^{v_{\infty}}(T),\infty)$ for some
$c^{v_{\infty}} \in (0,K)$. The function
 $v_{n+1}$ is the unique bounded solution (in the classical sense) of
\begin{equation}\label{eq:diff-eqn-1}
\begin{split}
& \A v_{n+1}(x,T)-  (\alpha+\lambda)\cdot v_{n+1} (x,T)+\lambda
\cdot (P v_{n})(x,T) - \frac{\partial}{\partial T} v_{n+1}(x,T)=0,
\quad x > c^{v_{n+1}}(T),
\\ & v_{n+1}(x,T) =(K-x), \quad x \leq c^{v_{n+1}}(T),
\end{split}
\end{equation}
and it
satisfies
\begin{equation}\label{eq:sm-fit1}
\frac{\partial}{\partial x}
v_{n+1}(x,T)\bigg|_{x=c^{v_{n+1}}(T)}=-1, \quad T>0.
\end{equation}
Moreover, $v_{\infty}$ is the unique bounded solution (in the classical sense) of
\begin{equation}\label{eq:diff-eqn-2}
\begin{split}
& \A v_{\infty}(x,T)-  (\alpha+\lambda)\cdot v_{\infty}
(x,T)+\lambda \cdot (P v_{\infty})(x,T) - \frac{\partial}{\partial
T} v_{\infty}(x,T)=0 \quad x >c^{v_{\infty}}(T),
\\ & v_{\infty}(x,T) =(K-x) \quad x \leq c^{v_\infty}(T),
\end{split}
\end{equation}
and it satisfies
\begin{equation}\label{eq:sm-fit2}
\frac{\partial}{\partial x}
v_{\infty}(x,T)\bigg|_{x=c^{v_{\infty}}(T)}=-1, \quad T>0.
\end{equation}
On the other hand,
\begin{equation}\label{eq:v-inf-marh-super}
\A v_{\infty}(x,T)-  (\alpha+\lambda)\cdot v_{\infty}
(x,T)+\lambda \cdot (P v_{\infty})(x,T) - \frac{\partial}{\partial
T} v_{\infty}(x,T) \leq 0 \quad x <c^{v_{\infty}}(T).
\end{equation}
\end{lemma}

\begin{proof}
The fact that $\C^{v_{n+1}}=(c^{v_{n+1}},\infty)$ and
$C^{v_{\infty}}=(c^{v_{\infty}},\infty)$ for some $c^{v_{n+1}} \in
(0,K)$ and $c^{v_{\infty}} \in (0,K)$ follows from
Lemma~\ref{lem:cont-reg} since the assumptions in
that lemma hold thanks to Corollaries~\ref{cor:increas},
\ref{cor:decreasing}; Remarks~\ref{rem:v-infty}
and~\ref{rem:sharp-upper-bound}; and Lemma~\ref{lem:fixed-po}.

The partial differential equations (\ref{eq:diff-eqn-1}),
(\ref{eq:diff-eqn-2}); and the inequality in \eqref{eq:v-inf-marh-super} are satisfied as a corollary of
Lemma~\ref{eq:lem-KS}; Corollaries~\ref{cor:increas} and
~\ref{cor:decreasing}, Remarks~\ref{rem:v-infty},
\ref{rem:sharp-upper-bound}; Lemmas~\ref{lem:fixed-po},
\ref{eq:v-infty-n-hold}.

Observe that since $v_{n}$ is convex
(Corollary~\ref{cor:increas}) and non-increasing
(Corollary~\ref{cor:decreasing}) with respect to its first
variable, $v_{n+1}$ $(=J v_{n})$
satisfies the smooth fit condition in (\ref{eq:sm-fit1}) as a
result of Lemma~\ref{lem:smooth-fit-cond}. The smooth fit
condition in (\ref{eq:sm-fit2}) holds for
$v_{\infty}$ as a result of
Lemma~\ref{lem:smooth-fit-cond} since
$v_{\infty}$ $(=Jv_{\infty})$
(Lemma~\ref{lem:fixed-po}) and $x \rightarrow
v_{\infty}(x,T)$ is non-increasing and
convex. 
\end{proof}

The next lemma will be used to verify the fact that $V=v_{\infty}$. The classical It\^{o}'s rule can not be applied to the process $t \rightarrow v_{\infty}(S_t,T-t)$ since the function $v_{\infty}$ may fail to be $C^{2,1}$ at $T \to c^{v_{\infty}}(T)$. As a result, the semi-martingale decomposition of the process $t \rightarrow v_{\infty}(S_t,T-t)$ may contain an extra term term due to the local time of the process $S$ at the free boundary.

\begin{lemma}\label{lem:ito}
Let $X=\{X_t; t \geq 0\}$ be a semi-martingale and $b: \R_+ \rightarrow \R$ be a continuous function of bounded variation. Let $F: \R \times \R_+ \rightarrow \R$ be a continuous function that is $C^{2,1}$ on $\bar{C}$ and $\bar{D}$ (it may not be necessarily $C^{1,2}$ across the boundary curve $b$), in which
\[
C \triangleq \{(x,t) \in \R \times \R_+: x<b(t)\}, \quad D \triangleq \{(x,t) \in \R \times \R_+: x>b(t)\}.
\]
That is, there exit two functions $F^{1}, F^{2}: \R \times \R_+ \rightarrow \R$, that $C^{2,1}$ on $\R \times \R_+$, and $F(x,t)=F^{1}(x,t)$ when $(x,t) \in C$ and $F(x,t)=F^{2}(x,t)$ when $(x,t) \in D$. Moreover, $F^{1}(b(t),t)=F^{2}(b(t),t)$.
Then the following generalization of It\^{o}'s formula holds:
\begin{equation}
\begin{split}
F(X_t,t)&=F(X_0,0)+\int_0^{t}\frac{1}{2}\left[F_t(X_{s-}+,s)+F_{t}(X_{s-}-,s)\right]ds
\\&+\frac{1}{2}\int_0^{t}\left[F_x(X_{s-}+,s)+F_x(X_{s-}-,s)\right]dX_s
+\frac{1}{2}\int_0^{t}1_{\{X_{s-} \neq b(s)\}}F_{xx}(X_{s-},s)d\left<X,X\right>_s^c
\\&+\sum_{0<s \leq t}\left\{F(X_s,s)-F(X_{s-},s)-\frac{1}{2}\Delta X_s \left[F_{x}(X_{s-}-,s)+F_{x}(X_{s-}+,s)\right]\right\}
\\&+ \frac{1}{2}\int_0^{t} \left[ F_x(X_{s-}+,s)-F_{x}(X_{s-}-,s)\right] 1_{\{X_{s-}=b(s)\}}dL_{t}^{b},
\end{split}
\end{equation}
where $L_t^b$ is the local time of the semi-martingale $X_t-b(t)$ at zero (see the definition on page 216 in \cite{Protter}).
\end{lemma}

Lemma~\ref{lem:ito} was stated in Theorem 2.1 of \cite{peskir05} for continuous semimartingales. The generalization for the case when the underlying process is not necessarily continuous is intuitively clear and just technical, but we will prove it in the Appendix for the sake of completeness. We are now ready to state the main results.

\begin{thm}\label{thm:main}
The value function $V$
is the unique bounded solution (in the classical sense) of the
integro-partial differential equation in (\ref{eq:diff-eqn-2}). Given
$(x,T)\in \R_+ \times \R_+ $ belongs to the optimal continuation
region if $x
>c^{v^{\infty}}(T)$.
Moreover, it satisfies the smooth fit condition at the optimal
stopping boundary, i.e.,
$\frac{\partial}{\partial x} V(x,T)\bigg|_{x=c^{v_{\infty}}(T)}=-1$, $T>0$.
\end{thm}
\begin{proof}
The proof is a corollary of the optional sampling theorem and the generalized It\^{o}'s formula given above. 
 Let $T \in (0,\infty)$ and define
\begin{equation}
\widetilde{M}_t=e^{-rt}v_{\infty}(S_t, T-t), \quad \text{and} \quad 
\widetilde{\tau}_x:=T \wedge \inf\{t \in [0,T]: S_t \leq c^{v_\infty}(T-t)\}.
\end{equation}
 It follows from (\ref{eq:diff-eqn-2}) and the classical It\^{o}'s
lemma that $\{\widetilde{M}_{t \wedge\widetilde{\tau}_x}\}_{0 \leq t \leq T}$ is a
bounded $\P$-martingale. Using the optional sampling theorem,
one obtains
\begin{equation}
v_{\infty}(x,T)=\widetilde{M}_0=\E \left\{\widetilde{M}_{\widetilde{\tau}_x}\right\}=\E\left\{e^{-r
\widetilde{\tau}_x}v_{\infty}(S_{\widetilde{\tau}_x},T-\widetilde{\tau}_x)\right\}=\E\left\{e^{-r
\widetilde{\tau}_x}(K-S_{\widetilde{\tau}_x})^+\right\} \leq V(x,T).
\end{equation}
In the rest of the proof we will show that $v_{\infty}(x,T) \geq V(x,T)$.
Since $v_{\infty}$ satisfies the smooth fit principle across the free boundary, when we apply the generalized It\^{o}'s formula to $v_{\infty}(S_t,T-t)$, the local time term drops. Thanks to \eqref{eq:diff-eqn-2} and \eqref {eq:v-inf-marh-super}, $v_{\infty}(S_t,T-t)$ is a
positive $\P$-supermartingale. Again, using the optional sampling
theorem, for any $\tau \in \widetilde{\S}_{0,T} $
\begin{equation}
v_{\infty}(x,T)=\widetilde{M}_0 \geq \E \left\{\widetilde{M}_{\tau}\right\}=\E\left\{e^{-r
\tau}v_{\infty}(S_{\tau},T-\tau)\right\} \geq \E\left\{e^{-r
\tau}(K-S_{\tau})^+\right\}.
\end{equation}
As a result $v_{\infty}(x,T) \geq V(x,T)$. 
\end{proof}

\begin{remark}
We have that
\begin{equation}\label{eq:cont-reg}
\C^{v_{\infty}}_T=\{x \in (0,\infty): v_{\infty}>(K-x)^+\}=(c^{v_{\infty}}(T),\infty).
\end{equation}
On the other hand, $v_{\infty}=K-x$ for $x \leq c^{v_{\infty}}$. Since $V=v_{\infty}$, by Theorem~\ref{thm:main}, it follows that  
\begin{equation}
\C^{V}_T=\{x \in (0,\infty):V>(K-x)^+\}=(c^{v_{\infty}}(T),\infty).
\end{equation}
\end{remark}

\setcounter{section}{0}%
\renewcommand{\thesection}{\Alph{section}}%
\section{Appendix}
\noindent \textbf{Proof of Lemma~\ref{lem:ito}.}
 As in \cite{peskir05} we will define 
$Z_t^1=X_t \wedge b(t)$, $ Z_t^2=X_t \vee b(t)$,
and observe that
\begin{equation}\label{eq:Fidentity}
F(X_t,t)=F^{1}(Z_t^1,t)+F^{2}(Z_t^2,t)-F(b(t),t).
\end{equation}
On the other hand, applying the Meyer-It\^{o} formula (see Theorem 70 in \cite{Protter}) to the semi-martingale $X_t-b(t)$, we obtain
\begin{equation}\label{eq:tanaka}
\begin{split}
|X_t-b(t)|&=|X_0-b(0)|+\int_0^t \text{sign} (X_{s-}-b(s))d(X_s-b(s))
\\&+2 \sum_{0 <s \leq t} \left[1_{\{X_{s-}>b(s)\}}(X_s-b(s))^{-}+1_{\{X_{s} \leq b(s) \}}(X_s-b(s))^{+}\right]+L_t^b.
\end{split}
\end{equation}  
Since 
$Z_t^1=\frac{1}{2}\left(X_t+b(t)-|X_t-b(t)|\right)$ and $Z_t^2=\frac{1}{2}\left(X_t+b(t)+|X_t-b(t)|\right)$,
using \eqref{eq:tanaka}, we get
\begin{equation}
\begin{split}
dZ_t^1&=\frac{1}{2}\left\{(1-\text{sign}(X_{t-}-b(t)))dX_t+(1+\text{sign}(X_{t-}-b(t)))db(t)-dL_t^b\right\}
\\&- \left[1_{\{X_{t-}>b(t)\}}(X_t-b(t))^{-}+1_{\{X_{t} \leq b(t) \}}(X_t-b(t))^{+}\right],
\end{split}
\end{equation}
\begin{equation}
\begin{split}
dZ_t^2&=\frac{1}{2}\left\{(1+\text{sign}(X_{t-}-b(t)))dX_t+(1+\text{sign}(X_{t-}-b(t)))db(t)-dL_t^b\right\}
\\&+\left[1_{\{X_{t-}>b(t)\}}(X_t-b(t))^{-}+1_{\{X_{t} \leq b(t) \}}(X_t-b(t))^{+}\right].
\end{split}
\end{equation}
It follows from the dynamics of $Z^i$, $i \in \{1,2\}$ that
\begin{equation}\label{eq:quad-var}
d\left<Z^i,Z^i\right>_t^c=\left(1_{\{X_{t-}<b(t)\}}+\frac{1}{4}1_{\{X_{t-}=b(t)\}}\right)d \left<X,X\right>_t^c
=1_{\{X_{t-}<b(t)\}}d \left<X,X\right>_t^c,
\end{equation}
where the second equality follows from the occupation density formula, see e.g. Corollary 1 on page 219 of \cite{Protter}. Applying the classical It\^{o}'s formula to $F^1(Z_t^1,t)$ and $F^{2}(Z_t^2,t)$ and using the dynamics
of $Z^1$ and $Z^2$, we get
\begin{equation}\label{eq:ito-F1}
\begin{split}
F^1(Z_t^1,t)&=F^1(Z_0^1,0)+\int_0^{t}F_t^1(Z_{s-}^1,s)ds+\int_0^{t}F_x^1(Z_{s-}^1,s)dZ_s^1+\frac{1}{2}
\int_0^{t}F_{xx}^1(s,Z_{s-}^1)d \left<Z^1,Z^1\right>^c_s
\\&+ \sum_{0  s \leq t}\left[F^1(Z_{s}^1,s)-F^1(Z_{s-}^1,s)-\Delta Z^1_s F^1_{x}(Z_{s-}^1,s)\right]
\\&=F^{1}(Z_0^1,0)+\int_0^{t}F_t^1(Z_{s-}^1,s)ds + \frac{1}{2}\int_0^{t} (1-\text{sign}(X_{s-}-b(s)))F_{x}^1(Z^1_{s-},s) dX_s
\\&+ \frac{1}{2}\int_0^{t} (1+\text{sign}(X_{s-}-b(s)))F_{x}^1(Z^1_{s-},s) db(s)
\\&-\sum_{0 <s \leq t} \left[1_{\{X_{s-}>b(s)\}}(X_s-b(s))^{-}+1_{\{X_{s} \leq b(s) \}}(X_s-b(s))^{+}\right] F_x^1(Z^1_{s-},s)
\\&-\frac{1}{2}\int_0^{t}F_{x}^1(Z^1_{s-},s)dL_t^b+\frac{1}{2}\int_0^{t} 1_{\{X_{s-}<b(s)\}}F_{xx}^1(Z^1_{s-},s)d\left<X^c,X^c\right>_s
\\& \sum_{0<s \leq t} \left[F^{1}(Z^1_s,s)-F(Z^1_{s-},s)-\Delta Z_{s}^1 F^1_{x}(Z^1_{s-},s)]\right], \quad \text{and}
\end{split}
\end{equation}
\begin{equation}\label{eq:ito-F2}
\begin{split}
F^2(Z_t^2,t)&=F^2(Z_0^2,0)+\int_0^{t}F_t^2(Z_{s-}^2,s)ds+\int_0^{t}F_x^2(Z_{s-}^2,s)dZ_s^2+\frac{1}{2}
\int_0^{t}F_{xx}^2(s,Z_{s-}^2)d \left<Z^2,Z^2\right>^c_s
\\&+ \sum_{0  s \leq t}\left[F^2(Z_{s}^2,s)-F^2(Z_{s-}^2,s)-\Delta Z^2_s F^1_{x}(Z_{s-}^2,s)\right]
\\&=F^{2}(Z_0^2,0)+\int_0^{t}F_t^2(Z_{s-}^2,s)ds + \frac{1}{2}\int_0^{t} (1+\text{sign}(X_{s-}-b(s)))F_{x}^2(Z^2_{s-},s) dX_s
\\&+ \frac{1}{2}\int_0^{t} (1-\text{sign}(X_{s-}-b(s)))F_{x}^2(Z^2_{s-},s) db(s)
\\&+\sum_{0 <s \leq t} \left[1_{\{X_{s-}>b(s)\}}(X_s-b(s))^{-}+1_{\{X_{s} \leq b(s) \}}(X_s-b(s))^{+}\right] F_x^2(Z^2_{s-},s)
\\&-\frac{1}{2}\int_0^{t}F_{x}^2(Z^2_{s-},s)dL_t^b+\frac{1}{2}\int_0^{t} 1_{\{X_{s-}<b(s)\}}F_{xx}^2(Z^2_{s-},s)d\left<X^c,X^c\right>_s
\\& +\sum_{0<s \leq t} \left[F^{2}(Z^2_s,s)-F(Z^2_{s-},s)-\Delta Z_{s}^2 F^2_{x}(Z^2_{s-},s)]\right].
\end{split}
\end{equation}
By splitting each term to their respective values on the sets $\{X_{s-}<b(s)\}$, $\{X_{s-}=b(s)\}$ and $\{X_{s-}>b(s)\}$, it can be seen that the following four equations are satisfied:
\begin{equation}\label{eq:peskir1}
F^{1}(Z_0^1,0)+F^{2}(Z_0^2,0)=F(X_0,0)+F(b(0),0),
\end{equation}
\begin{equation}\label{eq:peskir2}
\begin{split}
&\int_0^{t}F_{t}^{1}(Z_{s-},s)ds+\int_{0}^t F_{t}^2(Z_{s-}^2,s)ds
=\frac{1}{2}\int_{0}^{t}F_{t}(X_{s-}+,s)+F_{t}(X_{s-}-,s)ds
\\&+ \int_{0}^{t}\bigg[ F_{t}(b(s)+,s)1_{\{X_{s-}< b(s)\}}+\frac{1}{2}\left(F_{t}(b(s)-,s)+F_{t}(b(s)+,s)\right)1_{\{X_{s-}=b(s)\}}+F_{t}(b(s)-,s)1_{\{X_{s-}>b(s)\}}\bigg] ds,
\end{split}
\end{equation}
\begin{equation}\label{eq:peskir3}
\begin{split}
 \frac{1}{2}\int_0^{t} (1-\text{sign}(X_{s-}-b(s)))F_{x}^1(Z^1_{s-},s) dX_s &+\frac{1}{2}\int_0^{t} (1+\text{sign}(X_{s-}-b(s)))F_{x}^2(Z^2_{s-},s) dX_s 
 \\& =\frac{1}{2}\int_0^{t} \left[F_{x}(X_{s-}+,s)+F_{x}(X_{s-}-,s)\right]dX_s,
\end{split}
\end{equation}
\begin{equation}\label{eq:peskir4}
\begin{split}
& \frac{1}{2}\int_0^{t} (1+\text{sign}(X_{s-}-b(s)))F_{x}^1(Z^1_{s-},s) db(s)
+ \frac{1}{2}\int_0^{t} (1-\text{sign}(X_{s-}-b(s)))F_{x}^2(Z^2_{s-},s) db(s)=
\\&\int_0^{t} \bigg[ F_{x}(b(s)+,s)1_{\{X_{s-}<b(s)\}}+\frac{1}{2}\left[F_{x}(b(s)+,s)+F_{x}(b(s)-,s)\right]1_{\{X_{s-}=b(s)\}}+F_{x}(b(s)-,s)1_{\{X_{s-}>b(s)\}}\bigg]db(s).
\end{split}
\end{equation}
On the other hand, (3.15) of \cite{peskir05} still holds: 
\begin{equation}\label{eq:peskir5}
\begin{split}
F(b(t),t)&=F(b(0),0)+\int_{0}^{t}\bigg[F_{t}(b(s)+,s)1_{\{X_{s-}<b(s)\}} +\frac{1}{2}\left[F_{t}(b(s)-,s)+F_{t}(b(s)+,s)1_{\{X_{s-}=b(s)\}}\right]
\\&+F_{t}(b(s)-,s)1_{\{X_{s-}>b(s)\}}\bigg]ds
\\&\int_{0}^{t}\bigg[F_{x}(b(s)+,s)1_{\{X_{s-}<b(s)\}} +\frac{1}{2}\left[F_{x}(b(s)-,s)+F_{x}(b(s)+,s)1_{\{X_{s-}=b(s)\}}\right]
\\&+F_{x}(b(s)-,s)1_{\{X_{s-}>b(s)\}}\bigg]db(s),
\end{split}
\end{equation}
whose proof is carried out by using the uniqueness of finite measures on p-systems.

Let us analyze the jump terms in (\ref{eq:ito-F1}) and (\ref{eq:ito-F2}). We will denote
\begin{equation}
\begin{split}
A &:= -\left[1_{\{X_{s-}>b(s)\}}(X_s-b(s))^{-}+1_{\{X_{s} \leq b(s) \}}(X_s-b(s))^{+}\right] F_x^1(Z^1_{s-},s)
\\&+\left[F^{1}(Z^1_s,s)-F(Z^1_{s-},s)-\Delta Z_{s}^1 F^1_{x}(Z^1_{s-},s)]\right],
\end{split}
\end{equation}
\begin{equation}
\begin{split}
B &:=\left[1_{\{X_{s-}>b(s)\}}(X_s-b(s))^{-}+1_{\{X_{s} \leq b(s) \}}(X_s-b(s))^{+}\right] F_x^2(Z^2_{s-},s)
\\&+\sum_{0<s \leq t} \left[F^{2}(Z^2_s,s)-F(Z^2_{s-},s)-\Delta Z_{s}^2 F^2_{x}(Z^2_{s-},s)]\right].
\end{split}
\end{equation}
Depending on the whereabouts of $X_{s-}$ and $X_{s}$ with respect to the boundary curve $b$, $A$ and $B$ take four different values:
\begin{enumerate}
\item $X_{s-}>b(s)$ and $X_{t} \geq b(t)$. In this case 
\begin{equation}
A=0, \; B=F^{2}(X_{s},s)-F^{2}(X_{s-},s)-\Delta X_s F^{2}_x(X_{s-},s),
\end{equation}
and 
\begin{equation}\label{eq:effectofjump1}
A+B=F(X_s,s)-F(X_{s-},s)- \Delta X_{s} F_{x}(X_{s-}+,s).
\end{equation}
\item $X_{s-}>b(s)$ and $X_{s}<b(s)$. In this case
\begin{equation}
\begin{split}
A&=-(b(s)-X_{s})F_{x}^{1}(b(s),s)+F^{1}(X_{s},s)-F^{1}(b(s),s)-(X_s-b(s))F_x^{1}(b(s),s)
\\&=F^{1}(X_s,s)-F^{1}(b(s),s),
\end{split}
\end{equation} 
\ \begin{equation}
\begin{split}
B&=(b(s)-X_s)F^{2}_x(b(s),s)+F^{2}(b(s),s)-F^{2}(X_{s-},s)-(b(s)-X_{s-})F_{x}^2(X_{s-},s)
\\&=F^{2}(b(s),s)-F^{2}(X_{s-},s)- \Delta X_{s} F_{x}^2 (X_{s-},s), \quad \text{and}
\end{split}
\end{equation} 
\begin{equation}\label{eq:effectofjump2}
A+B=F(X_s,s)-F(X_{s-},s)-\Delta X_{s} F_{x}(X_{s}+,s).
\end{equation}
\item $X_{s-} \leq b(s)$ and $X_{s} \geq b(s)$. We have that
\begin{equation}
\begin{split}
A&=-(X_s-b(s))F_{x}^{1}(X_{s-},s)+F^{1}(b(s),s)-F^{1}(X_{s-},s)-(b(s)-X_{s-}) F^{1}_x(X_{s-},s)
\\&= F^{1}(b(s),s)-F^{1}(X_{s-},s)- \Delta X_{s} F^1(X_{s-},s),
\end{split}
\end{equation}
\begin{equation}
\begin{split}
B&=(X_s-b(s))F_x^2(b(s),s)+F^{2}(X_{s},s)-F^{2}(b(s),s)-(X_s-b(s))F_{x}^{2}(b(s),s)
\\&=F^{2}(X_s,s)-F^{2}(b(s),s).
\end{split}
\end{equation}
As a result
\begin{equation}\label{eq:effectofjump3}
A+B= F(X_{s},s)-F(X_{s-},s)- \Delta X_{s} F_{x}(X_{s-}-,s).
\end{equation}
\item $X_{s-} \leq b(s)$ and $X_{s}<b(s)$. Clearly,
\begin{equation}
A=F^{1}(X_s,s)-F^{1}(X_{s-},s)- \Delta X_{s} F_{x}^1(X_{s-},s) \quad \text{and} \quad B=0.
\end{equation}
As a result
\begin{equation}\label{eq:effectofjump4}
A+B= F(X_s,s)-F(X_{s-},s)- \Delta X_{s} F_x(X_{s-},s).
\end{equation}
\end{enumerate}
Now combining \eqref{eq:Fidentity}, \eqref{eq:quad-var}, \eqref{eq:ito-F1}, \eqref{eq:ito-F2}, \eqref{eq:peskir1},
\eqref{eq:peskir2}, \eqref{eq:peskir3}, \eqref{eq:peskir4}, \eqref{eq:peskir5},
\eqref{eq:effectofjump1},
\eqref{eq:effectofjump2},
\eqref{eq:effectofjump3},
\eqref{eq:effectofjump4}, we obtain
\begin{equation}\label{eq:almost-ito}
\begin{split}
F(X_t,t)&=F(X_0,0)+\frac{1}{2}\int_0^{t} \left[F_{t}(X_{s-}+,s)+F_{t}(X_{s-}-,s)\right] ds
\\&+\frac{1}{2}\int_{0}^{t}\left[F_{x}(X_{s-}+,s)+F_{x}(X_{s-}-,s)\right]dX_s+\frac{1}{2}\int_0^{t}1_{\{X_{s-} \leq b(s)\}}F_{xx}(s,X_{s-})d \left<X,X\right>^c_{s-}
\\& +\sum_{0<s \leq t}\bigg[F(X_s,s)-F(X_{s-},s)-\Delta X_s F_{x}(X_{s-}-,s)1_{\{X_{s-} \leq b(s)\}}- \Delta X_{s} F_{x}(s,X_{s-}+)1_{\{X_{s-}>b(s)\}}\bigg]
\\&+\frac{1}{2}\int_0^{t} \left[F_{x}^2(Z^2_{s-},s)-F^1(Z^1_{s-},s)\right]dL_{t}^b.
\end{split}
\end{equation}
The last term on the right-hand-side of (\ref{eq:almost-ito}) can be written as
\begin{equation}
\frac{1}{2}\int_0^{t} \left[F_{x}^2(Z^2_{s-},s)-F^1(Z^1_{s-},s)\right]dL_{t}^b= \frac{1}{2} \int_0^{t} \bigg[F_{x}(X_{s-}+,s)-F_{x}(X_{s-}-,s)\bigg]1_{\{X_{s-}=b(s)\}}dL_t^b,
\end{equation}
using Theorem 69 of \cite{Protter}. On the other hand, the jump term in  \eqref{eq:almost-ito} can be written as 
\begin{equation}
\begin{split}
\sum_{0<s \leq t}&\bigg[F(X_s,s)-F(X_{s-},s)-\Delta X_s F_{x}(X_{s-}-,s)1_{\{X_{s-} \leq b(s)\}}- \Delta X_{s} F_{x}(s,X_{s-}+)1_{\{X_{s-}>b(s)\}}\bigg]
\\& =\sum_{0<s \leq t}\bigg[F(X_s,s)-F(X_{s-},s)-\frac{1}{2}\Delta X_s \left[F_x(X_{s-}-,s)+F_{x}(X_{s-}+,s)\right] \bigg].
\end{split}
\end{equation}
This completes the proof. \hfill $\square$

{\small
\bibliographystyle{plain}
\bibliography{references}
}

\end{document}